\documentclass[12pt]{amsart}
\usepackage{amsmath,amsthm}

  \def\ee{{\bf e}}

\def\pmatrix{\left(\begin{matrix}}
\def\endpmatrix{\end{matrix}\right)}

\def\R{{\mathbb R}}
\def\P{{\mathbb P}}
\def\Z{{\mathbb Z}}
\def\C{{\mathbb C}}

\def\Q{{\mathbb Q}}

\def\cal{\mathcal}

\def\de{\delta}

\def\p{\partial}

\def\t{\theta}
\def\s{\sigma}
\def\T{\Theta}
\def\e{\varepsilon}

\def\b{\beta}

\def\A{{\mathcal A}}

\def\H{{\mathcal H}}

\def\tt#1#2{{\t\left[\begin{matrix}{#1}\\ {#2}\end{matrix}\right]}}

\def\nZ#1{\left(\frac{1}{#1}\Z/\Z\right)^g}
\theoremstyle{plain}
\newtheorem{thm}{Theorem}
\newtheorem{lm}[thm]{Lemma}
\newtheorem{prop}[thm]{Proposition}

\theoremstyle{definition}

\newtheorem{rem}[thm]{Remark}

\title{Theta functions of arbitrary order and their derivatives}
\author{Samuel Grushevsky\and Riccardo Salvati Manni}
\address{Mathematics Department, Princeton University, Fine Hall,
Washington Road, Princeton, NJ 08544, USA}
\email{sam@math.princeton.edu}
\address{Dipartimento di Matematica, Universit\`a di Roma ``La
Sapienza'', Piazzale Aldo Moro, 2, I-00185 Roma, Italy}
\email{salvati@mat.uniroma1.it}
\thanks{First author partially supported by the NSF Mathematical
Sciences Postdoctoral Research Fellowship}

\begin{document}
\begin{abstract}
In this paper we establish the relationships between theta functions
of arbitrary order and their derivatives. We generalize our previous
work \cite{paper1} and prove that for any $n>1$ the map sending an
abelian variety to the set of Gauss images of its points of order
$2n$ is an embedding into an appropriate Grassmannian (note that for
$n=1$ we only got generic injectivity in \cite{paper1}). We further
discuss the generalizations of Jacobi's derivative formula for any
dimension and any order.
\end{abstract}
\maketitle

\section{Introduction and definitions}
The study of theta functions of abelian varieties is a very
classical subject that goes back to Jacobi, Riemann, Weierstrass,
Fr\"obenius, Poincar\'e and many others. A purely algebraic modern
treatment of the subject started with Weil \cite{weil}. In the
1960s Igusa \cite{iginj} and Mumford \cite{mumford} proved the
fundamental theorem relating the values of theta functions at zero
to injective maps from some modular varieties into the projective
space, among other results. For a detailed history of the problem
up till 1980 we refer to Igusa's survey \cite{igreview}.

\smallskip
More precisely, let $\H_g$ be the Siegel upper half-space --- the
set of symmetric $g\times g$ complex matrices $\tau$ with
positive-definite imaginary part. For any $z\in\C^g$,
$\tau\in\H_g$ we define the theta function with characteristics
$\e,\de\in\R^g$ to be
$$
\tt\e\de(\tau,z)=\sum\limits_{n\in\Z^g}\ee\left(\frac 12
(n+\e)^t \tau(n+\e)+(n+\e)^t(z+\de)\right),
$$
where $\ee(t):=\exp(2\pi i t)$, and $A^t$ denotes the transpose of
a matrix $A$. The restriction of a theta function to $z=0$ is
called the associated theta constant.

If in the formula above we take $\e\in\nZ{n},$ $\de=0,$ and
replace $\tau$ and $z$ by $n\tau$ and $nz$, the resulting theta
function is called the theta function with characteristics of
order $n$. These theta functions form a basis for the space of
sections of $n\T$ --- the $n$-th power of the symmetric line
bundle inducing the principal polarization on the abelian variety
with period matrix $\tau$. In this case the associated theta
constants are modular forms with respect to a certain subgroup of
${\rm Sp}(2g,\Z)$. Let us define this.

The symplectic group ${\rm Sp}(2g,\Z)$ acts on $\H_g$. Let us
write an element $\gamma\in{\rm Sp}(2g,\Z)$ as $\gamma=\pmatrix
a&b\\ c&d\endpmatrix$ with $a,b,c,d$ being $g\times g$ integer
matrices. Then the action is given by
$$
  \pmatrix a&b\\ c&d\endpmatrix \tau:=(a\tau+b)(c\tau+d)^{-1},
$$
and the quotient is the moduli space of principally polarized
abelian varieties $\A_g={\rm Sp}(2g,\Z)\backslash \H_g$. Let $\rho:{\rm
GL}(g,\C)\to \operatorname{End} V$ be an irreducible rational
representation with the highest weight $(k_1,k_2,\dots,k_g)$,
$k_1\geq k_2 \geq\dots\geq k_g$; then we call $k_g$ the weight of
$\rho$. A representation $\rho_0$ is called reduced if its weight
is equal to zero. Let us fix an integer $r$; we are interested in
pairs $\rho=(\rho_0,r)$, with $\rho_0$ reduced. We call $r$ the
weight of $\rho$ and use the notation
$$
  \rho(A)=\rho_0(A)\det A^{r/2}\ .
$$
For a finite index subgroup
$\Gamma\subset\operatorname{Sp}(2g,\Z)$ a multiplier system of
weight $r/2$ is a map $v:\Gamma\to \C^*$, such that the map
$$
  \sigma\mapsto v(\sigma)\det(C\tau+D)^{r/2}
$$
satisfies the cocycle condition for every $\sigma\in\Gamma$ and
$\tau\in\H_g$ (note that the function $\det(C\tau+D)$ possesses a
square root). Clearly a multiplier system of integral weight is a
character. A map $f:\H_g\to V$ is called a $\rho$- or $V$-valued
modular form, or simply a vector-valued modular form, if the
choice of $\rho$ is clear, with multiplier $v$, with respect to a
subgroup of finite index $\Gamma\subset\operatorname{Sp}(2g,\Z)$
if the transformation formula
\begin{equation}\label{transform}
  f(\sigma\circ\tau)=v(\sigma)\rho(C\tau+D)f(\tau)
\end{equation}
is satisfied for any $\sigma$ in $\Gamma$ and any $\tau$ in
$\H_g$, and, for $g=1$, if additionally $f$ is holomorphic at all
cusps of $\Gamma\backslash\H_1$.

Let us now define the level subgroups of the symplectic group to
be
$$
  \Gamma_g(n):=\left\lbrace\gamma=\pmatrix a&b\\ c&d\endpmatrix
  \in\Gamma_g\, |\, \gamma\equiv\pmatrix 1&0\\
  0&1\endpmatrix\ {\rm mod}\ n\right\rbrace,
$$
$$
  \Gamma_g(n,2n):=\left\lbrace\gamma\in\Gamma_g(n)\, |\, {\rm diag}
  (a^tb)\equiv{\rm diag}(c^td)\equiv0\ {\rm mod}\ 2n\right\rbrace.
$$
We denote the corresponding level covers of $\A_g$ by
$\A_g(n):=\Gamma_g(n)\backslash\H_g$ and
$\A_g(n,2n):=\Gamma_g(n,2n)\backslash\H_g$, respectively. It is
known that theta constants of order $n$ are modular forms of weight
$1/2$ (and with $\rho_0={\rm Id}$), for a suitable multiplier $v_n$,
with respect to the group $\Gamma_g(n,2n)$.

\smallskip
One of the main results proved by Igusa in \cite{iginj} and
Mumford in \cite{mumford} is that the map
$$
  Th_{n}:\A_g(n,2n)\to\P^{n^g-1}
$$
sending a point to the set of values of all theta constants of a
given order  $n$ for any $n\ge 4$,
$$
  Th_n(\tau):=\left\lbrace \tt a0(n\tau,0)\right\rbrace_{ {\rm all}
  \ a\in\nZ{n}},
$$
defines an embedding of the level moduli space.

Recently in \cite{paper1} we considered, in the case of both
characteristics $\e,\de$ of a theta function being half-integral
(which is equivalent to the order 4 case, see \cite{igusa}) the
map sending a point in $\H_g$ to the $g\times 2^{g-1}(2^g-1)$
matrix of non-trivial gradients
$$
  \Phi_4:\tau\to\left({\rm grad}_{z=0}\tt\e\de(\tau,z)\right)_{ {\rm
  all\ odd}\ \e,\de\in\nZ{2}}
$$
and showed that $\Phi_4$ induces a generically injective and
immersive map of $\A_g(4,8)$ to the Grassmannian variety
$G(g,2^{g-1}(2^g-1))$ of $g$-planes in $\C^{2^{g-1}(2^g-1)}$. Here
``odd'' means that as a function of $z$ the theta function
$\tt\e\de(\tau,z)$ is odd, which is equivalent to the scalar
product $2(\e,\de)$ being zero in $\nZ{2}$.

Passing to Pl\"ucker's coordinates to embed the Grassmannian
variety into a projective space, the image of the map $\Phi_4$ in
the projective space produces some well-known modular forms, the
so-called Jacobian determinants of theta functions. These are
obtained as follows: for any set of $g$ odd characteristics
$[\e_1,\de_1],\dots, [\e_g,\de_g]$ we define their jacobian
determinant to be
$$
  \begin{matrix}
    D([\e_1,\de_1],\dots[\e_g,\de_g])(\tau):=\qquad\qquad\hfill\\
    \qquad\hfill \pi^{-g} grad\,\,\tt{\e_{1}}{\de_{1}}(\tau,0)\wedge
    grad\,\,\tt{\e_{2}}{\de_{2 }}(\tau,0)\wedge \dots \wedge
    grad\,\,\tt{\e_{g}}{\de_{g}}(\tau, 0).
  \end{matrix}
$$

The Jacobian determinants were also extensively studied in the
nineteenth century, with special emphasis on their modular
properties and relationship with theta constants. The first result
in that direction was the famous Jacobi's derivative formula
$$
  \tt{\frac 12}{\frac 12}'(\tau,0)=-\pi\tt00(\tau,0)
  \tt{\frac 12}0(\tau,0)\tt0{\frac 12}(\tau,0),
$$
which is the expression of the only non-zero Jacobian determinant
for $g=1$ in terms of even theta constants.

Generalizations of this formula were stated by Rosenhain in the case
of genus two \cite{rosenhain},  and mostly proved by Weber and
Fr\"obenius for genera up to four \cite{frobenius}. It seems that
Riemann has also worked on this problem, and some generalizations
can be found in \cite{riemann}. Thomae \cite{thomae} then
generalized the formula to the case of hyperelliptic curves of any
genus, but the problem of completely generalizing Jacobi's
derivative formula to arbitrary abelian varieties remained open.

Recently in \cite{paper2} we found different generalizations of
Jacobi's derivative formula to higher genus, involving second
order derivatives of theta functions at zero.

\medskip
The aim of this paper is to present a general framework for
deriving the generalizations of the results of
\cite{paper1},\cite{paper2} to arbitrary level. We will consider
the map
$$
  \Phi_{4n}:\A_g(4n,8n)\to G(g,n^{2g})
$$
given by
$$
  \Phi_{4n}(\tau):=\left\lbrace {\rm grad}_{z=0}\tt a0(4n\tau,z)
  \right\rbrace_{{\rm all}\ a\in\nZ{4n}}
$$
and relate it to the theta constant maps $Th_{2n}$ and $Th_{4n}$.
Note that the map $\Phi_{4n}$ is well-defined because gradients of
theta functions of order $4n$ are vector-valued modular forms with
respect to $\Gamma_g(4n,8n)$ for the representation ${\rm
std}\otimes\det^{1/2}$ (i.e. of weight $1/2$ and with
$\rho_0(A)=A$). Notice that unlike the $n=1$ case, here for
convenience we include the gradients of all theta functions
irrespective of their parity, though of course since $\tt
a0(n,z)+\tt{-a}0(n,z)$ is even, there will be many identical
columns in the $g\times n^{2g}$ matrix, which is the image
$\Phi_{4n}(\tau)$.

\smallskip
We will show that $\Phi_{4n}$ is an embedding for all $n>1$
(recall that in \cite{paper1} we considered the case of $n=1$ and
were only able to prove generic injectivity), and will also obtain
generalizations of Jacobi's derivative formula for theta functions
of arbitrary level. We think that similar results can also be
obtained for other levels not divisible by 4, but dealing with
those makes some computations much more technically involved, as
working with theta functions of non-integral level is harder, and
we will not treat such computations here.

\smallskip
We will work with theta functions of orders $2n$ and $4n$ (in
\cite{paper1} and \cite{paper2} we worked with $n=1$). To try to
avoid confusion, we will adhere to the following notations: Greek
letters will stand for characteristics $\e\in\nZ{2}$, which will
play a special role, Latin letters at the end of the alphabet will
be for vectors $z\in\C^g$, and Latin letters at the beginning of
the alphabet will be for characteristics $a\in\nZ{m}$ for some
even order $m$, or sometimes for $a\in(\Q/\Z)^g$ for complete
generality.

\section{Addition theorem for theta functions}
We work with a $g$-dimensional principally polarized abelian
variety $X=V/\Lambda$ with period matrix $\tau$ and the
polarization bundle $\T$. We denote by $X[2]$ the points of order
two $X$, i.e. points $p\in X$ such that $2p=0\in X$. For $x=\tau
\e+\de$ in $X[2]$ the shifted bundle $t^*_x\T$ is still a
symmetric line bundle. The theta function $\tt\e\de(\tau,z)$ is,
up to a multiplicative constant, the unique section of $t^*_x\T$.
Note, however, that $\tt\e\de(\tau,2z)$ is a section of $\T^4$ due
to the presence of the lower characteristic. In general a basis of
$H^0(X,\T^n)$ is given by the $n^g$ theta functions $\tt a
0(n\tau,nz)$ with $a\in\nZ{n}$. We now recall the formula in
\cite{igusa} at the top of p.~50:
\begin{equation}
\label{1}
  \tt a b( \tau, z+\tau c +d)=\ee\left(-\frac{1}{2}c^{t}(\tau c+
 z+d+b)\right)\, \tt {a+c}{b+d}(\tau,z).
\end{equation}
We will also need a slight generalization of  the formula at the
bottom of p.~171 in \cite{igusa}, relating theta functions of
order twice larger and theta functions with a lower
characteristic:

\begin{lm}
\label{know1} For all $\tau\in\H_g,\ z\in\C^g,$ $a\in\R^g$ and
$\b\in\nZ{2}$ we have
$$
  \tt a\b(\tau,2z)=\sum\limits_{\e\in\nZ{2}}\ee(\b^t(2\e+a))
  \tt{\e+\frac{a}2}0(4\tau,4z).
$$
\end{lm}

One of the basic relations among theta functions is Riemann's
bilinear addition theorem, which essentially relates theta
functions at $\tau$ and $2\tau$ or, if the characteristics are
chosen appropriately, theta functions of order $n$ and $2n$. We
will need to use it in two forms. The first form is the following

\begin{prop}[specialization of Theorem 2, p. 139 in \cite{igusa}]
\label{know} For all $\tau\in\H_g,\ z,w\in\C^g,$ $a,b\in \R^g,$
and $\e\in\nZ{2}$ the following holds:
$$
  \tt a\e(4n\tau,4nz)\tt b\e(4n\tau,4nw)=
$$
$$
  =\frac{1}{2^g}\sum\limits_\s\ee(-2a^t\s)\,
  \tt{a+b}{\s+\e} (2n\tau,2n(z+w))\tt{a-b}{\s}(2n\tau,2n(z-w)).
$$
\end{prop}

We will also need another form of this addition theorem, which in
some sense is the converse, expressing one term in the
right-hand-side of the above as a combination of terms in the
left-hand-side.

\begin{prop}[a generalization of \cite{igusa}, Corollary, p. 141]
\label{addnew}
For all $\tau\in\H_g,\ z,w\in\C^g,\ a,b\in \R^g,$
and $\gamma,\s\in\nZ{2}$ the following holds
$$
  \tt a{\gamma+\s}(2n\tau,2nz)\, \tt
  b\gamma(2n\tau,2nw)=\sum\limits_\e\ee
  ((a+b+2\e)^t\gamma)\cdot
$$
$$
  \cdot \tt{\e+\frac{a+b}2}\s(4n\tau,2n(z+w))\,
  \tt{\e+\frac{a-b}2}\s(4n\tau,2n(z-w)).
$$
\end{prop}

\begin{proof}
This formula differs from the one in the previous proposition in
that we are trying to pass to double argument rather than half the
argument. We first apply formula (\ref{1}) to the left-hand-side
and then use proposition \ref{know}. Afterwards we use the formula
in lemma \ref{know1}.
$$
  \tt a{\gamma+\s}(2n\tau,2nz)\,
  \tt b\gamma(2n\tau,2nw)=
$$
$$
  \tt a0(2n\tau,2nz+\gamma+\s)\,
  \tt b0(2n\tau,2nw+\gamma)=
$$
$$
  =\frac {1}{2^g}\sum\limits_\mu\ee (-2a^t\mu)\,
  \tt{a+b}\mu\left(n\tau,n (z+w)+\gamma+\frac\s 2\right)\cdot$$
$$
  \tt{a-b}\mu\left(n\tau,n (z-w)+\frac\s 2\right)
$$
$$
  =\frac{1}{2^g}\sum\limits_{\mu,\e,\de}\ee\left(-2a^t\mu+
 (2\e+a+b)^t\mu+ (2\de+a-b)^t\mu\right)\cdot
$$
$$
  \cdot\tt{\e+\frac{a+b}2}0(4n\tau,2n(z+w)+2\gamma+\s)\,
  \tt{\de+\frac{a-b}2}0(4n\tau,2n(z-w)+\s).
$$
When we take the sum over $\mu$ in this formula, this is just
taking the sum $\sum_\mu \ee(2(\e+\de)^t\mu)$, which is zero
unless $\e=\de$ and is equal to $2^g$ if $\e=\de$. Thus summing
over $\mu$ extracts $2^g$ times the $\e=\de$ terms of the above
sum, and we end up with
$$
  \sum\limits_\e\, \tt{\e+\frac{a+b}2}0(4n\tau,2n(z+w)+2\gamma+\s)\,
  \tt{\e+\frac{a-b}2}0(4n\tau,2n(z-w)+\s)
$$
$$
  =\sum\limits_\e\ee((a+b+2\e)^t\gamma)\,
  \tt{\e+\frac{a+b}2}\s(4n\tau,2n(z+w))\,
  \tt{\e+\frac{a-b}2}\s(4n\tau,2n(z-w))
$$
\end{proof}

We end this section by recalling that as $a$ varies in $\nZ{n}$,
by formula (\ref{1}) the values of the theta functions $\tt a 0(
n\tau, nz)$ at $0$ are related to the values of the single theta
function $\tt 0 0(n\tau,nz)$ at different points of order $n$ on the
abelian variety.

\section{Injectivity of the gradient maps}
In this section we follow, generalize and further advance the
framework of establishing the relationships between  gradients of
theta functions and derivatives of theta constants that we have
developed in \cite{paper1} and \cite{paper2}. We then use the
general relationships between the maps $\Phi$ and $Th$ to show
that the image of $Th$ can be obtained uniquely from the image of
$\Phi$, thus eventually proving injectivity of $\Phi_{4n}$ for
$n>1$. The improvement over the $n=1$ case, where we could only
get generic injectivity, is due to the fact that we can now
preclude the massive vanishing of theta constants that plagued our
computations in \cite{paper1}; we are also aided by the knowledge
that $Th_{2n}$ is an embedding for $n>1$, while it is still only a
conjecture that $Th_2$ is injective.

For simplicity, we denote by $\p_i\theta$ the derivative of
$\theta$ with respect to $z_i$, evaluated at $z=0$. Similarly to
\cite{paper1} and \cite{paper2}, let us then define the $g\times
g$ matrices
$$
  {\bf C}^{ab}:=\left(2\p_i\tt a0(4n\tau)\p_j\tt b0(4n\tau)+
  2\p_j\tt a0(4n\tau)\p_i\tt b0(4n\tau)\right)_{
  {\rm all}\ i,j}
$$
for $a,b\in (\Q/\Z)^g$. We mainly shall use $\bf C$ with both
indices $a,b\in\nZ{4n}$. Note that the $\bf C$'s that we used in
\cite{paper1} and \cite{paper2} are essentially the case $n=1$ of
the above, but here we used different indices for $\bf C$'s, since
we are using a different basis for theta functions of a given
order. Let us also define the $g\times g$ matrices
$$
  {\bf A}^{cd}_{\e}:=\left(\p_i\p_j \tt{c}\e(2n\tau)\,\tt{d}
  \e(2n\tau)-\tt{c}\e(2n\tau)\p_i\p_j\tt{d}\e(2n\tau)\right)_{ {\rm
  all}\ i,j}
$$
for $a,b\in(\Q/\Z)^g$ and $\e\in\nZ{2}$. Similarly to $\bf C$, the
$\bf A$'s we used in our previous works correspond to the case
$n=1$ of the definition we are now using, with some further
restrictions on $a$ and $b$ .

Note that $\bf A$ and $\bf C$ are vector-valued modular forms with
respect to $\Gamma_g(4n,8n)$ and the representation
$$
  \rho= Sym^2 ({\rm std})\otimes\det.
$$

\begin{thm}
\label{AC}
The matrices $\bf A$ and $\bf C$ can be expressed in terms of
each other as follows:

a) ${\bf C}^{ab}=\frac{1}{2^g}\sum\limits_\s\ee(-2a^t\s){\bf A}^{
a+b,a-b}_{\s}.$

b) ${\bf A}^{ab}_{\de}=2 \sum\limits_\e\ee((a+b+2\e)^t\de)
{\bf C}^{\e+\frac{a+b}{2},\e+\frac{a-b}{2}}.$
\end{thm}

\begin{proof}
Indeed, to get part a) let us take the derivative
$\p_{z_i}\p_{w_j}+\p_{z_j}\p_{w_i}$ of the formula in proposition
\ref{know}, and then evaluate at $z=w=0$. Differentiating the
left-hand-side is easy. On the right-hand-side we notice that the
terms where each factor is differentiated once will cancel because
of the minus sign for $w$ in the argument of the second theta
function. Thus we arrive at
$$
  2\p_i\tt a\e(4n\tau)\p_j\tt b\e(4n\tau)+
  2\p_j\tt a\e(4n\tau)\p_i\tt b\e(4n\tau)
$$
$$
  =\frac{1}{2^g}\sum\limits_\s\ee(-2a^t\s)\left(\p_i\p_j
  \tt{a+b}{\s+\e}(2n\tau)\,\tt{a-b}\s(2n\tau)\right.
$$
$$
  -\left.\tt{a+b}{\s+\e}(2n\tau)\p_i\p_j\tt{a-b}\s(2n\tau)\right),
$$
which, when written in terms of $\bf A$ and $\bf C$, gives us part
a) of the theorem.

For the proof of part b) let us take the derivative
$\p_{z_i}\p_{z_j}-\p_{w_i}\p_{w_j}|_{z=w=0}$ of the formula in
proposition \ref{addnew}. Differentiating the left-hand-side is
easy; on the right-hand-side we notice that the terms that do not
cancel are the ones where each of the factors is differentiated
once, and thus we end up with
$$
  \p_i\p_j\tt a\gamma(2n\tau)\,\tt b\de(2n\tau)- \tt
  a\gamma(2n\tau)\p_i\p_j\tt b\de(2n\tau)=
$$
$$
  2\sum\limits_\e\ee((a+b+2\e)^t\gamma)\cdot
  \left( \p_i\tt{\e+\frac{a+b}2}{\gamma+\de}(4n\tau)\p_j
  \tt{\e+\frac{a-b}2}{\gamma+\de}(4n\tau)\right)+
$$
$$
  2\sum\limits_\e\ee((a+b+2\e)^t\gamma)\cdot
  \left(\p_j\tt{\e+\frac{a+b}2}{\gamma+\de}(4n\tau)\p_i
  \tt{\e+\frac{a-b}2}{\gamma+\de}(4n\tau)\right),
$$
which in terms of $\bf A$ and $\bf C$ is exactly part b) of the
theorem.
\end{proof}

In the following we will only use this theorem for the case when
the indices of ${\bf C}$ lie in $\nZ{4n}$ and the upper indices of
${\bf A}$ lie in $\nZ{4n}$ with the extra condition that the
indices of the corresponding ${\bf C}$ appearing in part a) of
theorem \ref{AC} are all in $\nZ{4n}$, i.e. with the condition
that $a+b\in \nZ{2n}$.

\smallskip
We will now proceed to show the injectivity of the gradient theta
map at all levels --- this is done similarly to the computations
in \cite{paper1} while taking advantage of the more general $\bf
A$ and $\bf C$, so we now streamline the argument.

\begin{lm}
The following identity holds:
$$
  {\bf A}^{ab}_\e\,\tt c\e(2n\tau)+{\bf A}^{bc}_\e\,
  \tt a\e(2n\tau)+{\bf A}^{ca}_\e\,\tt b\e(2n\tau)=0.
$$
\end{lm}

\begin{proof}
This is a trivial computation with all the six terms canceling
pairwise.
\end{proof}

We observe that the above lemma in particular holds for $a, b,
c\in\nZ{4n}$ with $a+b\in \nZ{2n}$ and $a+c\in\nZ{2n}$

In an improvement over the $n=1$ case, where we had trouble
proving non-degeneracy, we can now prove

\begin{lm}
For any $n>1$ the rank of the $(2n)^g\times (\frac{g(g+1)}{2}+1)$
matrix with columns
$$
  \left(\tt {a + \frac{\delta} {2n}}\e(2n\tau),\p_i\p_j\tt  {a
  +\frac{\delta}{2n}}\e(2n\tau)\right)_{{\rm all}\ a\in \nZ{2n},
  \ {\rm all}\ (i,j)}
$$
for any fixed $\e$ and $\delta$ is maximal, for all $g\ge 1$.
\end{lm}

\begin{proof}
Lemma 11, p. 188 in \cite{igusa} proves this result for
$\delta=\e=0$ and for any order $m$ divisible by $4$. The proof
given there clearly works for any even $m=2n\ge 4$ as well. Using
formula (\ref{1}), we can then obtain a proof of the lemma by
evaluating the theta functions $\tt a0(2n\tau, 2nz)$ at the point
$z=\tau \frac{\delta}{2n}+ \frac{\e}{2n}$.
\end{proof}

The reason why the $n=1$ case would not work for the lemma above
is that all even theta functions vanish at odd points. We would
also like to remark that this result is closely related to the
injectivity of certain higher order embeddings of abelian
varieties --- obtained by using theta functions, not their
derivatives --- which were studied in \cite{bs}.

\smallskip
Now similarly to proposition 12 in \cite{paper1} we can
reconstruct the (projectivized) values of theta constants from the
knowledge of ${\bf A}$'s and thus, by theorem \ref{AC}, from the
${\bf C}$'s, i.e. from $\Phi_{4n}(\tau)$.

\begin{prop} \label{uniq2n}
The value of $\Phi_{4n}(\tau)$ uniquely determines for any fixed
$\gamma,\de\in\nZ{2}$ the projective point
$$
  \left\lbrace\tt {a +\frac{\delta}{2n}}\gamma(2n\tau)\right
  \rbrace_{{\rm all}\ a\in\nZ{2n}}.
$$
\end{prop}

For $\gamma=0$ this point is simply the value $Th_{2n}(\tau)$.
Since we know that $Th_{2n}$ is an embedding of $\A_g(2n,4n)$ for
$n>1$, this means that $\Phi_{4n}(\tau)$ determines the class of
$\tau$ in $\A_g(2n,4n)$ uniquely. Since the cover $\A_g(4n,8n)\to
\A_g(2n,4n)$ is finite, it follows immediately that the map
$\Phi_{4n}$ on $\A_g(4n,8n)$ is at most finite-to-one. We would
now like to show that $\Phi_{4n}$ is in fact injective by showing
that from the knowledge of $\Phi_{4n}(\tau)$ we can determine
uniquely the class of $\tau$ in $\A_g(4n,8n)$ and not only in
$\A_g(2n,4n)$. The first step in this direction is the following

\begin{thm}
For any fixed $\sigma,\de\in\nZ{2}$ and fixed $a, b\in\nZ{2n}$
with $a+b-\frac{\delta}{n}\in \nZ{n}$, the value of
$\Phi_{4n}(\tau)$ uniquely determines the projective point
$$
  \left\lbrace\tt{a }{\gamma+\s}(n\tau)\tt{ b}{\gamma}(n\tau)
  \right\rbrace_{{\rm all}\ \gamma}.
$$
\label{uniqprod}
\end{thm}

\begin{proof}
Indeed, let us use the addition formula from proposition
\ref{addnew} with $2n$ and $4n$ replaced by $n$ and $2n$. Then in
the right-hand-side we will have a linear combination of terms
appearing in proposition \ref{uniq2n}, which are uniquely
determined by $\Phi_{4n}(\tau)$, while in the left-hand-side we
will be getting products of two theta functions at $n\tau$ of the
kind described.
\end{proof}

The problem we had in \cite{paper1} in trying to prove injectivity
was due in large part to the possibility of many theta constants
vanishing simultaneously, so that we were unable to determine
certain signs uniquely. For $n>1$ we can deal with this.

\begin{lm}
For all $n>1$ and for any fixed $\gamma,\s,\delta\in\nZ{2}$ there
always exist some $a,b\in\nZ{2n}$ with $a+b-\frac{\delta}{n}\in
\nZ{n}$, such that
$$
  \tt a{\gamma+\sigma}(n\tau)\tt b\gamma(n\tau)\neq 0.
$$
\end{lm}

\begin{proof}
First note that for any fixed $\gamma+\sigma$ there is at least
one among $\tt a{\gamma+\sigma}(n\tau)$ that does not vanish:
indeed, these are the values of all theta functions of order $2n$
at the point $\gamma+\sigma$. Thus let us pick some $a$ such that
$\tt a{\gamma+\sigma}(n\tau)\ne 0$.

Similarly  let us consider theta functions of order $n$, $\tt
c0(n\tau,nz)$ for $c\in\nZ{n}$. Among these there is at least one
not vanishing at $z= \tau (-a+ \frac{\delta}{n})+
\frac{\gamma}{n}$; let us choose such a $c$. We then finally set
$b:=c-a+\frac{\delta}{n}$, and by formula \ref{1} this implies
that $\tt {b}\gamma (n\tau,0)\neq 0.$
\end{proof}

Now we are ready to prove the main result.

\begin{thm}
The map $\Phi_{4n}$ is injective on $\A_g(4n,8n)$ for all $n>1$
and all $g\geq 2$.
\end{thm}

\begin{proof}
Recall that $\Gamma_g(2n,4n)/\Gamma_g(4n,8n)$ acts on theta
constants of order $4n$ by multiplying them by $\pm 1$, depending
on characteristics. Thus to finish reconstructing $Th_{4n}(\tau)$
from $\Phi_{4n}(\tau)$ (and thus also knowing $Th_{2n}(\tau)$) we
need to deal with the ``projectivization'' happening in theorem
\ref{uniqprod}, to recover the necessary signs. By the formulas on
page 171 of \cite{igusa} (see also section 2 of this paper),
instead of considering the theta constants $\tt c0(4n\tau)$ with
$c\in\nZ{4n}$, we can consider the theta constants $\tt
c\mu(n\tau)$ with $c\in\nZ{2n}$.

Indeed, suppose that $\Phi_{4n}(\tau)=\Phi_{4n}(\tau')$. The
previous lemma states that for fixed $\s,\delta,\gamma$ we can
always find $a$ and $b$ with $a+b-\frac{\delta}{n}\in \nZ{n}$ such
that
$$
  \tt a{\gamma+\sigma}(n\tau)\tt b\gamma(n\tau)\neq 0.
$$
Since such products are projectively unique by theorem
\ref{uniqprod}, we have
$$
  \tt a{\gamma+\sigma}(n\tau)\tt b\gamma(n\tau)=t_{\s , \,\delta}\,
  \tt a{\gamma+\sigma}(n\tau')\tt b\gamma(n\tau')
$$
for some (unique, since both sides are non-zero --- this is
crucial!) constant $t_{\s , \,\delta}$ independent of $\gamma$.

Squaring the above formula we get
$$
 t_{\s , \,\delta}^2=t_ {0 , \,0} ^2
$$

We claim then that the map
$$
  X:\left(\frac{1}{2}\Z/\Z\right)^{2g}\to\pm 1,\qquad
  X(\s,\,\delta):=t_{\s , \,\delta}/t_ {0 , \,0}
$$
is a group morphism. In fact for fixed $\s$ we can always find
$a,b\in\nZ{n}$ with
$$
  \tt a{ \sigma}(n\tau)\tt b 0(n\tau)\neq 0
$$
and for fixed $\delta$ we can always find $ b_1$ with
$b_1-\frac{\delta}{n}\in \nZ{n}$ satisfying
$$
  \tt b{ 0}(n\tau)\tt{ b_1} 0(n\tau)\neq 0.
$$
This then shows that $X(\s,\,\delta)= X(\s,\, 0) X(0,\,\delta).$
Similarly one proves that $X(\s+\rho,\, 0) =X(\s,\, 0)X(\rho,\, 0)
$ and $X(0,\,\delta+\e)=X(0,\,\delta)X(0,\,\e),$ and thus we see
that $X$ is indeed a morphism.

\smallskip
\noindent To show that $Th_{4n}(\tau)=Th_{4n}(\tau')$, we need to
show that $X$ is identically equal to +1. Since $X$ is a morphism,
we only need to check that a basis gets mapped to $+1$. If this is
not the case, then we have some $X(\s,\,\delta)=-1$. Then we can
find an element $M\in\Gamma_g(2n,4n)/\Gamma_g(4n,8n)$ (in fact
such an element can be found in $\Gamma_g(4n)/\Gamma_g(4n,8n)$),
the action of which on theta constants of level $4n$ would change
precisely the appropriate signs --- the argument for $n>1$ is
identical to the one given in \cite{paper1} for $n=1$.

Thus we know that
$$
  \tt a{\gamma+\sigma}(n\tau)\tt b\gamma(n\tau)=
  \tt a{\gamma+\sigma}(nM\tau')\tt b\gamma(nM\tau'),
$$
from which it follows that $Th_{4n}(\tau)=Th_{4n}(M\tau')$ --- by
fixing some $b,\gamma$ such that $\tt b\gamma(n\tau)\ne0$ and
varying $a$ and $\gamma$, so that we get $\tt
a{\gamma+\sigma}(n\tau)/\tt a{\gamma+\sigma}(nM\tau')={\rm const}$
independent of $a$ and $\sigma$. Since $Th_{4n}$ is injective, it
means that $\tau$ and $M\tau'$ represent the same point in
$\A_g(4n,8n)$. This then implies that $\Phi_{4n}(\tau')=
\Phi_{4n}(\tau)=\Phi_{4n}(M\tau')$. However, there cannot be an
$M\not\in\Gamma_g(4n,8n)$ such that its action does not change the
image under the map $\Phi_{4n}$ (see \cite{paper1}). Hence we must
have $M\in\Gamma_g(4n,8n)$, so $\tau=\tau'$ in $\A_g(4n,8n)$, and
thus the  injectivity of $\Phi_{4n}$ is proved.
\end{proof}

\begin{rem}
We observe that the assumption $n>1$ has been used to prove that
$X(\cdot\, ,\,\cdot)$ is a homomorphism. In fact, for $n=1$ we
could not show that $X$ is indeed defined, as we did not have the
non-vanishing results and thus some of $t_{\s,\de}$ could be
undefined if many theta constants vanished.
\end{rem}

\begin{rem}
The injectivity of $\Phi_{4n}$ on the tangent spaces follows from
lemma 17 in \cite{paper1} and the result in \cite{SM1}.
\end{rem}

\section{Generalized Jacobi's derivative formulas}
In the same spirit as above, the results of \cite{paper2} can be
generalized to higher level. The relationship between $\bf A$ and
$\bf C$ provides us with a way to express vector-valued modular
forms constructed using theta constants and their
$\tau$-derivatives (which, by the heat equation, are the same as
the second $z$-derivatives) in terms of the gradients of theta
functions. These can be used to deduce relations among scalar
modular forms involving Jacobian determinants of theta functions.
In fact both formulas from \cite{paper2} can be generalized to
higher level. Below we give the appropriate version of Theorem 5
from that paper.

We recall the matrix differential operator
$$
  \cal D:= \left(\begin{array}{rrrr}
  \,\frac{\p}{\p\tau_{11}}&\frac{1}{2}\frac{\p}
  {\p\tau_{12}}&\dots&\frac{1}{2}\frac{\p}{\p\tau_{1 g}}\\
  \frac{1}{2}\frac{\p}{\p \tau_{21}}&\frac{\p}{\p
  \tau_{22}}&\dots&\frac{1}{2}\frac{\p}{\p\tau_{2 g}}\\
  \dots&\dots&\dots&\dots\\
  \frac{1}{2}\frac{\p}{\p \tau_{g 1}}& \dots&\dots& \,\
  \frac{\p}{\p\tau_{g g}}\end{array}\right).
$$
Then we have

\begin{thm}
For any $a\in\nZ{2n},\de\in\nZ{2}$ the following holds:
\begin{equation}
\label{Eq}
  \begin{matrix}
    {\rm const}\ \left(\tt 0\de(2n\tau)\right)^{2g}\det
    \left(\cal D(\tt a\de(2n\tau)/\tt 0\de(2n\tau))\right)\\
    =\!\!\!\!\!\sum\limits_{\e_1,\ldots,\e_g\in\nZ{2}}\!\!\!\!\!\ee(2\de^t
    (\e_1+\ldots+\e_g)) D([a/2+\e_1,0],\dots [a/2+\e_g,0]
    )^2(4n\tau)
  \end{matrix}
\end{equation}
for some computable constant ${\rm c}$.
\end{thm}

\begin{proof}
This follows by linear algebra arguments from the expression of
$\bf A$ in terms of $\bf C$ and applying the Binet's formula to
the matrix ${\bf C}^{aa}$, which has rank one, being equal to the
product of a vector and a covector. The proof is the same as in
\cite{paper2}.
\end{proof}

\end{document}